\documentclass{proc-l}

\newtheorem{theorem}{Theorem}[section]
\newtheorem{lemma}[theorem]{Lemma}

\theoremstyle{definition}
\newtheorem{definition}[theorem]{Definition}

\newtheorem{corollary}[theorem]{Corollary}

\theoremstyle{remark}
\newtheorem{remark}[theorem]{Remark}

\numberwithin{equation}{section}

\begin{document}

\title[On Dixmier and Connes-Dixmier traces]{On the distinction between the classes of Dixmier and Connes-Dixmier traces}


%    author one information
\author{Fedor Sukochev}
\address{School of Mathematics and Statistics, University of New South Wales, Sydney, 2052, Australia.}
\curraddr{}
\email{f.sukochev@unsw.edu.au}
\thanks{}

%    author two information
\author{Alexandr Usachev}
\address{}
\curraddr{}
\email{a.usachev@unsw.edu.au}
\thanks{}

%    author three information
\author{Dmitriy Zanin}
\address{}
\curraddr{}
\email{d.zanin@unsw.edu.au}
\thanks{Research supported by the Australian Research Council}

%    \subjclass is required.

\date{}

\dedicatory{}

%    "Communicated by" -- provide editor's name; required.
\commby{Varghese Mathai}

%    Abstract is required.
\begin{abstract}
 In the present paper we prove that the classes of Dixmier and Connes-Dixmier traces differ even on the
Dixmier ideal $\mathcal M_{1,\infty}$. We construct a Marcinkiewicz space $\mathcal M_\psi$ and a positive operator $T\in \mathcal M_\psi$ which is Connes-Dixmier measurable but which is not Dixmier measurable.
\end{abstract}

\maketitle

\section{Introduction and preliminares}
In~\cite{D} J.~Dixmier proved that there exists a non-normal trace (a Dixmier trace) on the non-commutative Marcinkiewicz spaces $\mathcal M_\psi$ 
for every $\psi$ such that 
\begin{equation}\label{cond_lim}
\lim_{t\to\infty} \frac{\psi(2t)}{\psi(t)}=1. 
\end{equation}
 In \cite{C} A.~Connes introduced a subclass of Dixmier traces, later termed in~\cite{LSS} Connes-Dixmier traces. 
In this paper, we investigate the relationship between these two classes and show that they differ even on the classical Dixmier ideal 
$\mathcal M_{1,\infty}$. 
Furthermore, we prove that there is a Marcinkiewicz ideal $\mathcal M_\psi$, with $\psi$ satisfying~\eqref{cond_lim} such that these two classes of traces generate distinct sets of measurable elements (see~\cite[IV.2.$\beta$.Definition 7]{C} and Definitions~\ref{Dmeas} and~\ref{CDmeas} below).

\subsection{Generalized limits}
Let $l_\infty$ 
be the Banach space of all bounded sequences $x=(x_0,x_1,\ldots)$ with the norm
$$\|x\|_{l_\infty}:=\sup_{n\ge0} |x_n|.$$ 

A normalized positive linear functional on $l_\infty$ which equals the ordinary limit on convergent sequences is called a generalized limit.
For every $n\in \mathbb N$ we define a dilation operator $\sigma_n : l_\infty \to l_\infty$ as follows
$$\sigma_n (x_0, x_1, \ldots) = \left(\underbrace{x_0,\ldots ,x_0}_{n},\underbrace{x_1,\ldots ,x_1}_{n},\ldots \right).$$
If a generalized limit $\omega$ on $l_\infty$ satisfies the condition
$$\omega(\sigma_n x)=\omega(x)$$ 
for every $x\in l_\infty$ and any $n\in \mathbb N$, then $\omega$ is called a dilation invariant generalized limit.

Let $L_\infty=L_\infty(0,\infty)$ be the space of all real-valued
bounded Lebesgue measurable functions on $(0,\infty)$ equipped with the norm
$$\|x\|_{L_\infty}:=\mathop{\rm esssup}\limits_{t>0} |x(t)|.$$

A normalized positive linear functional on $L_\infty$ which equals the ordinary limit on convergent (at infinity) sequences is called a generalized limit.
For every $x\in L_\infty$ and for any generalized limit $\gamma$ on $L_\infty$ the following inequalities hold
$$\liminf_{t\to\infty} x(t)\le \gamma(x)\le \limsup_{t\to\infty} x(t).$$
By Hahn-Banach extension theorem, for every $x\in L_\infty$ there exist generalized limits $\gamma_1$ and $\gamma_2$ such that
\begin{equation}\label{gl}
\gamma_1(x)= \limsup_{t\to\infty} x(t) , \ \ \gamma_2(x)= \liminf_{t\to\infty} x(t).
\end{equation}

We define a dilation operator $\sigma_s : L_\infty \to L_\infty$ as follows
$$(\sigma_s x)(t) = x(t/s), \quad s>0.$$
A generalized limit $\omega$ on $L_\infty$ is said to be dilation invariant if 
$$\omega(\sigma_s x)=\omega(x)$$ 
for every $x\in L_\infty$ and any $s>0$.

Let $\pi$ be the isometric embedding $\pi:l_{\infty}\to L_\infty$
given by 
$$\{x_n\}_{n=0}^\infty\stackrel{\pi}{\mapsto}\sum_{n=0}^\infty x_n \chi_{(n,n+1]}.$$

The following natural way to generate dilation invariant generalized limits was suggested in~\cite[Section IV, 2$\beta$]{C}.  A.~Connes 
observed that for any generalised limit $\gamma$ on $L_\infty$ 
a functional $\omega:=\gamma\circ M\circ \pi$ is a dilation invariant generalized limit on $l_\infty$.
Here, the bounded operator $M:L_\infty\to L_\infty$ is given by the formula
$$(Mx)(t):=\frac1{\log t} \int_1^t x(s)\,\frac{ds}s.$$

Throughout the paper we denote by $\log t$ the natural logarithm and by $\log_2 t$ the logarithm with base 2.

\subsection{Marcinkiewicz spaces}
Let $B(H)$ be an algebra of all bounded linear operators on a separable Hilbert space
$H$ equipped with the uniform norm and let ${\rm Tr}$ be the standart trace.

For every operator $T\in B(H)$ a generalized singular value function $\mu(T)$ is defined by the formula
$$\mu(t, T)=\inf \{\|Tp\| : \ p \ \text{is a projection in} \ B(H) \ \text{with} \ {\rm Tr}(1-p)\le t\}.$$
For a compact operator $T$, it can be proven that $\mu(k, T)$ is the $k$-th largest eigenvalue of an operator $|T|$, $k \ge 0$.

Since $B(H)$ is an atomic von Neumann algebra and traces of all atoms equal to 1, it follows that 
$\mu(T)$ is a step function and $\mu(T)=\pi(\mu(k,T))$ for every $T\in B(H)$.

Let $\Omega$ denote the set of all concave functions $\psi: [0,\infty) \to [0,\infty)$ such that 
$\lim_{t\to0+}\psi(t)=0$ and $\lim_{t\to\infty}\psi(t)=\infty$.

Let $\psi \in \Omega$. Consider the Banach ideal $(\mathcal M_\psi,\|\cdot\|_{\mathcal M_\psi})$ of compact operators in $B(H)$ given by (see e.g.~\cite{CS,KSS,LSS})
$$\mathcal M_\psi:= \left\{T: \|T\|_{\mathcal M_\psi}:= \sup_{n\ge 0}  \frac{1}{\psi(1+n)} \sum_{k=0}^n\mu(k,T) < \infty\right\}.$$

For $f\in L_\infty$ we set 
$$a(t,f):=\frac1{\psi(t)}\int_0^t f^*(s)ds,$$
where $f^*$ denotes the decreasing rearrangement of the function $|f|$ that is
$$f^*(t):=\inf\{s\geq0:\ mes(\{|f|>s\})\leq t\}.$$

We define the Marcinkiewicz function space $M_\psi$
of real-valued measurable functions $f$ on $(0,\infty)$ by setting
$$\|f\|_{M_\psi}:=\sup_{t>0}a(t,f)<\infty.$$

For a compact operator we have $T\in \mathcal M_\psi$ if and only if $\mu(T) \in M_\psi$.

In the case when $\psi(t)=\log(1+t)$ the space $\mathcal M_\psi$ is a well-known Dixmier ideal $\mathcal M_{1,\infty}$.

\subsection{Singular traces on general Marcinkiewicz spaces}

For an arbitrary dilation invariant generalized limit $\omega$ on $l_\infty$ 
the weight
$$
{\rm Tr}_\omega (T):= \omega \left( \left\lbrace \frac{1}{\log(2+n)} \sum_{k=0}^n\mu(k,T)\right\rbrace _{n=0}^\infty\right), \quad 0\le T\in \mathcal M_{1,\infty},
$$
extends to a non-normal trace (a Dixmier trace) on $\mathcal M_{1,\infty}$~\cite{D,C,CS}. We denote the set of all Dixmier traces by $\mathcal D$.

The subclass $\mathcal{C}\subset \mathcal{D}$ of all Dixmier traces ${\rm Tr}_{\omega}$ defined by $\omega=\gamma\circ M\circ \pi$ 
was termed Connes-Dixmier traces in~\cite{LSS}. A priori, $\mathcal{C}\subseteq \mathcal{D}$ 
and the question about precise relationship between these two classes arises naturally.
 Recently the distinction between $\mathcal{C}$ and $\mathcal{D}$ was studied by A.~Pietsch 
in terms of density characters (see~\cite{P1}-\cite{P3}). 
For the discussion of various classes of singular traces we refer to~\cite{CPS2, CS, LS}.

The first main result of the present paper (Theorem~\ref{CD_in_D} below) shows that the inclusion $\mathcal{C} \subset \mathcal{D}$ is proper.
Our approach is completely different from that of A.~Pietsch and the proof provided here is much shorter.

It has become traditional to reduce various problems about Dixmier traces to its commutative analogues.

For every dilation invariant generalized limit $\omega$ on $L_\infty$  one can define a commutative analogue of Dixmier trace 
(a Dixmier functional on $M_{1,\infty}$) as follows
\begin{equation}\label{Df}
\tau_\omega(f)=\omega(a(t,f)), \ 0\le f\in M_{1,\infty} 
\end{equation}
and extend it to $M_{1,\infty}$ by linearity.

It was shown in~\cite{DPSS, KSS} that, for a general Marcinkiewicz space $\mathcal M_\psi$, the following conditions are equivalent
\begin{enumerate}
\item[(i)] The space $\mathcal M_\psi$ admits non-trivial Dixmier traces.
\item[(ii)] The function $\psi \in \Omega$ satisfies the following condition
\begin{equation}\label{cond}
 \liminf_{t\to\infty} \frac{\psi(2t)}{\psi(t)}=1.
\end{equation} 
\item[(iii)] There exists a dilation invariant limit $\omega$ on $l_\infty$ such that
\begin{equation}\label{comp1}
\omega\left(\frac{\psi(2n)}{\psi(n)}\right)=1. 
\end{equation}

\end{enumerate}

It was also proven in~\cite[Proposition 9, Theorem 11]{KSS} that for $\psi \in \Omega$ satisfying~\eqref{cond},
the weight
$$
{\rm Tr}_\omega (T):= \omega \left( \left\lbrace \frac{1}{\psi(1+n)} \sum_{k=0}^n\mu(k,T)\right\rbrace _{n=0}^\infty\right), \quad 0\le T\in \mathcal M_\psi,
$$
extends to a Dixmier trace on $\mathcal M_\psi$ if and only if 
a dilation invariant generalized limit $\omega$ on $l_\infty$ satisfies~\eqref{comp1}.

Similarly to the definition of Connes-Dixmier traces on $\mathcal M_{1,\infty}$, for every $\psi \in \Omega$ satisfing~\eqref{cond} and 
any dilation invariant limit $\omega=\gamma\circ M\circ \pi$ on $l_\infty$ satisfying~\eqref{comp1}
we can define a Connes-Dixmier trace ${\rm Tr}_\omega$ on $\mathcal M_\psi$.

Similarly to~\eqref{Df}, we define Dixmier and Connes-Dixmier functionals $\tau_\omega$ 
for every dilation invariant generalized limit $\omega$ on $L_\infty$ satisfying
\begin{equation}\label{comp2}
\omega\left(\frac{\psi(2t)}{\psi(t)}\right)=1. 
\end{equation}

\begin{remark}\label{duality}
By~\cite[Theorem 14, Corollary 15]{SZ} and [8, Theorem 11] we know that for every Dixmier trace ${\rm Tr}_{\omega_1}$ on $\mathcal M_\psi$
($\omega_1$ is a dilation invariant generalized limit on $l_\infty$)
there exists a Dixmier functional $\tau_{\omega_2}$ on $M_\psi$ ($\omega_2$ is a dilation invariant generalized limit on $L_\infty$) such that
$${\rm Tr}_{\omega_1}(T)= \tau_{\omega_2}(\mu(T)), \quad 0\le T\in \mathcal M_\psi. $$
The converse implication also holds.
\end{remark}

The following Lemma was borrowed from~\cite[Proposition 9]{KSS}.
\begin{lemma}\label{lem1}
Let $\psi \in \Omega$ satisfy~\eqref{cond} and let $\omega$ be a dilation invariant generalized limit on $L_\infty$ satisfying~\eqref{comp2}.
For every $f\in M_\psi$, we have
\begin{equation}
\omega\left(\frac{tf^*(t)}{\psi(t)}\right)=0. 
\end{equation} 
\end{lemma}

\begin{proof}
 Since $\omega$ is a dilation invariant generalized limit,
\begin{align*}
\omega\left(\frac{1}{\psi(2t)}\int_0^{2t}f^*(s)\,ds\right)&=\omega\left(\frac{1}{\psi(t)}\int_0^{t}f^*(s)\,ds\right)\\
&= \omega\left(\frac{\psi(2t)}{\psi(t)}\frac{1}{\psi(2t)}\int_0^{t}f^*(s)\,ds\right).
\end{align*}

 Since $\omega$ satisfies~\eqref{comp2}, it follows from~\cite[Proposition 4]{KSS} that
$$\omega\left(\frac{1}{\psi(2t)}\int_0^{2t}f^*(s)\,ds\right)= \omega\left(\frac{1}{\psi(2t)}\int_0^{t}f^*(s)\,ds\right).$$

Hence,
$$\omega\left(\frac{1}{\psi(2t)}\int_t^{2t}f^*(s)\,ds\right)=0$$
and, furthermore,
$$\omega\left(\frac{2tf^*(2t)}{\psi(2t)}\right)=0.$$
Again, applying dilation invariance of $\omega$, we have
$$\omega\left(\frac{tf^*(t)}{\psi(t)}\right)=0.$$ 
\end{proof}

\subsection{Measurability}
The following definitions were motivated by A.~Connes~\cite[IV.2.$\beta$.Definition 7]{C} (see also~\cite[Definition 3.2]{LSS})
in the case when $\psi(t)=\log(1+t)$.
\begin{definition}\label{Dmeas}
Let $\psi \in \Omega$ satisfy~\eqref{cond}.
An operator $T\in \mathcal M_\psi$ is called Dixmier measurable if ${\rm Tr}_\omega(T)$ takes the same value for all ${\rm Tr}_\omega\in \mathcal{D}$.
\end{definition}
% Similarly we define the Connes-Dixmier measurability of an operator $T\in \mathcal M_\psi$ (see~\cite[Definition 5.9]{LSS}).

\begin{definition}\label{CDmeas}
Let $\psi \in \Omega$ satisfy~\eqref{cond}.
An operator $T\in \mathcal M_\psi$ is called Connes-Dixmier measurable if ${\rm Tr}_\omega(T)$ takes the same value for all ${\rm Tr}_\omega\in \mathcal{C}$.
\end{definition}

\begin{theorem}[Corollary 3.9 from~\cite{LSS}]
If $\psi\in \Omega$ satisfying~\eqref{cond} is such that
\begin{equation}\label{taub}
t\cdot\frac{d}{dt}\log\left(\psi(e^t)\right)<C 
\end{equation}
 for some $C>0$ and for all $t>0$, then for a positive operator $T\in \mathcal M_\psi$
the following statements are equivalent:

(i) $T$ is Dixmier measurable;

(ii) $T$ is Connes-Dixmier measurable;

(iii) There exists 
\begin{equation}\label{l0}
\lim_{n\to\infty} \frac{1}{\psi(n+1)} \sum_{k=0}^n\mu(k,T). 
\end{equation}
\end{theorem}

It is easy to check that the function $\psi(t)=\log(1+t)$ satisfies the condition~\eqref{taub}. 
Notwithstanding the difference between the sets of Dixmier and Connes-Dixmier traces, 
a positive operator $T\in \mathcal M_{1,\infty}$ is Connes-Dixmier measurable if and only if it is Dixmier measurable.

This result naturally raises the question, whether for an arbitrary function $\psi \in \Omega$ satisfying~\eqref{cond} the Connes-Dixmier measurability 
is equivalent to Dixmier measurability on the cone of all positive elements from $\mathcal M_\psi$.
Our second main result (Theorem~\ref{4} below) shows that the answer is (surpisingly) negative.

An example of the function $\psi\in \Omega$ satisfying~\eqref{cond} but failing the equivalence (i) $\Leftrightarrow$ (iii) was constructed in~\cite[Theorem 4.6]{DPSSS1}. 
However, if $\psi \in \Omega$ satisfies~\eqref{cond_lim}, 
then Theorem~\ref{lemma1} below shows that the equivalence (i) $\Leftrightarrow$ (iii) holds independently of the condition~\eqref{taub}.

\section{The classes of Dixmier and Connes-Dixmier traces are distinct}

Denote by $\mathcal M_\psi^0$ the separable part of the space $\mathcal M_\psi$, 
that is the closure in $\mathcal M_\psi$ of the set of all finite dimensional operators from $B(H)$.
The following Lemma was proven in~\cite[Theorems 2.8 and 5.12]{LSS} (see also~\cite[Theorems 7.3 and 7.4]{CS}).
\begin{lemma}\label{dist_lem}
If $\psi\in \Omega$ satisfies~\eqref{cond_lim}, then
\begin{equation}\label{d1}
{\rm dist}(T,\mathcal M_\psi^0)= \sup_{{\rm Tr}_\omega \in \mathcal{D}} {\rm Tr}_\omega(T) , \qquad 0\le T\in \mathcal M_\psi.
\end{equation}
If $\psi$ satisfies~\eqref{taub}, then there exists $c>1$ such that 
\begin{equation}\label{d2}
\sup_{{\rm Tr}_\omega \in \mathcal{C}} {\rm Tr}_\omega(T) \le {\rm dist}(T,\mathcal M_\psi^0) \le c \cdot \sup_{{\rm Tr}_\omega \in \mathcal{C}} {\rm Tr}_\omega(T), 
\qquad 0\le T\in \mathcal M_\psi.
\end{equation} 
\end{lemma}

In view of the difference between~\eqref{d1} and~\eqref{d2}, the following question arises naturally: 
"Is the constant $c$ in~\eqref{d2} necessarily strictly greater than 1?" 
The following theorem shows that the inclusion $\mathcal{C} \subset \mathcal{D}$ is proper and answers this question in the affirmative.

\begin{theorem}\label{CD_in_D}
There exists a positive operator $T_0\in \mathcal M_{1,\infty}$
such that
$$\sup_{{\rm Tr}_\omega \in \mathcal{D}} {\rm Tr}_\omega(T_0) > \sup_{{\rm Tr}_\omega \in \mathcal{C}} {\rm Tr}_\omega(T_0).$$
\end{theorem}
\begin{proof} Let $T_0$ be such that
$$\mu(T_0)=\sup_{k\geq0}2^{k-2^k}\chi_{[0,2^{2^k})}.$$ 
We set $f_0=\mu(T_0)$. 
By Lemma~\ref{dist_lem}, we have
\begin{align*}
\sup_{{\rm Tr}_\omega \in \mathcal{D}} {\rm Tr}_\omega(T_0)&= {\rm dist}(T,\mathcal M_\psi^0)\\
&=\limsup_{t\to\infty} a(t,f_0).
\end{align*}

By~\eqref{gl}, we have 
$$\sup_{{\rm Tr}_\omega \in \mathcal{C}} {\rm Tr}_\omega(T_0) = \limsup_{t\to\infty} (Ma(\cdot,f_0))(t).$$

So, it is sufficient to prove that
$$\limsup_{t\to\infty} a(t,f_0)> \limsup_{t\to\infty} (Ma(\cdot,f_0))(t).$$

Clearly, $f_0=f_0^*.$ For every $2^{2^n}\le t< 2^{2^{n+1}}$ we have
\begin{equation}\label{a1}
\begin{aligned}
a(t,f_0)&=\frac1{\log(1+t)}\left(\int_0^{2^{2^n}}f_0(s)ds+(t-2^{2^n})f_0(t)\right)\\
&=\frac1{\log(1+t)}\left(\sum_{k=1}^{n}\int_{2^{2^{k-1}}}^{2^{2^k}}2^{k-2^k}ds+tf_0(t)+O(1)\right)\\
&=\frac{2^{n+1}}{\log t}+\frac{tf_0^*(t)}{\log(1+t)}+o(1).
\end{aligned} 
\end{equation}

It is easy to check that $f_0\in M_{1,\infty}$ and, hence, $T_0\in \mathcal M_{1,\infty}$.

By Lemma~\ref{lem1}, $\gamma \circ M\left(\frac{tf^*_0(t)}{\log(1+t)}\right)=0$ for every generalized limit $\gamma$ on $L_\infty$
and, appealing to~\eqref{gl}, we conclude 
\begin{equation}\label{eq1}
\lim_{t\to\infty} M\left(\frac{sf^*_0(s)}{\log(1+s)}\right)(t)=0, \ \ \text{for every} \ f\in M_{1,\infty}. 
\end{equation}

Define the function $x\in L_\infty$ by setting
$$x(t):=\sum_{n=0}^\infty\frac{2^n}{\log t}\chi_{[2^{2^n},2^{2^{n+1}})}(t) \quad t>0.$$
Hence, we obtain from~\eqref{a1} and~\eqref{eq1}
$$\limsup_{t\to\infty} (Ma(\cdot,f_0))(t)=2\limsup_{t\to\infty} (Mx)(t).$$

For ${2^{2^n}}\le t <{2^{2^{n+1}}}$, we have
$$(M x)(t)= \frac1{\log t} \left(\sum_{k=0}^{n-1} 2^k\int_{2^{2^k}}^{2^{2^{k+1}}} \frac{d\log s}{\log s} + 2^n\int_{2^{2^n}}^t \frac{d\log s}{\log s} +O(1)\right).$$
Since
$$\int \frac{d\log s}{\log s} =\frac{\log_2(\log_2 s)-\log_2(\log_2 e)}{\log_2 e}+C,$$
it follows that
\begin{align*}(M x)(t)&= \frac{\log2}{\log t} \left(\sum_{k=0}^{n-1} 2^k + 2^n(\log_2(\log_2 t) -n)\right)+O(1)\\
&= \frac{2^n\log2}{\log t} \left(1 + \log_2(\log_2 t) -n\right)+o(1).\\
\end{align*}

 The function 
$$g:t \to \frac{2^n\log2}{\log t} \left(1 + \log_2(\log_2 t) -n\right), \quad t\in [2^{2^n}, 2^{2^{n+1}}) $$ 
has extrema at $$t_n=2^{2^{\frac1{\log2}-1 +n}}\in [2^{2^n}, 2^{2^{n+1}}), \ n\in \mathbb N.$$
We have $g(t_n)=\frac2{e \log2}$ for every $n\in \mathbb N$. 
Since  $g(2^{2^n})=1$ for every $n\in \mathbb N$ and since $g$ is continuous on $(1,\infty)$, 
it follows that $\limsup_{t\to\infty} g(t)=\frac2{e \log2}$
and
$$\limsup_{t\to\infty} (Ma(\cdot,f_0))(t)=2\limsup_{t\to\infty} (Mx)(t)=\frac4{e \log2}.$$

By the definition we have $f_0(2^{2^n})=2^{n+1-2^{n+1}}$ and so, from~\eqref{a1} we obtain 
\begin{align*}
\limsup_{t\to\infty} a(t,f_0) &\ge \limsup_{n\to\infty}a(2^{2^n},f_0)\\
&=\frac1{\log2} \cdot\limsup_{n\to\infty}\left(\frac{2^{n+1}}{\log_2 2^{2^n}} +\frac{2^{2^n}f_0(2^{2^n})}{\log_2 2^{2^n}}\right)\\
&=\frac1{\log2} \cdot\limsup_{n\to\infty}\left(2+\frac{2^{2^n}2^{n+1-2^{n+1}}}{2^n}\right)\\
&=\frac2{\log2}\\ 
&>\frac4{e \log2}\\ 
&=\limsup_{t\to\infty} (Ma(\cdot,f_0))(t).
\end{align*}
\end{proof}

\section{The classes of Dixmier and Connes-Dixmier measurable elements are distinct}

The following Lemma is taken from~\cite{SUZ} (see Theorem 18 or~\cite[Theorem 6.1.3]{Sed}).
\begin{lemma} Let $x \in L_\infty$ such that $x\circ exp$ is uniformly continuous.
The equality $\omega(x)=A$ holds for every dilation invariant generalized limit $\omega$ on $L_\infty$ if and only if
$$\lim_{t\to\infty}\frac1{\log t}\int_1^t x(\alpha s) \frac{ds}s =A$$
uniformly in $\alpha\ge1$.
\end{lemma}

\begin{corollary}\label{add0} Let $\psi \in \Omega$ satisfies~\eqref{cond}. Let $f\in M_\psi$ and let $A$ be a real number.
The equality $\tau_\omega(f)=A$ holds for every Dixmier functional $\tau_\omega$ if and only if
$$\lim_{t\to\infty} \frac1{\log t} \int_1^t a(\alpha s,f) \frac{ds}s =A $$
uniformly in $\alpha\ge1$.
\end{corollary}
\begin{proof}
The mapping $t\to a(e^t,f)$ is uniformly continuous since
\begin{align*}
 \left|\frac{d}{dt}(a(e^t,f))\right|&=\left|\frac{d}{dt}\left(\frac1{\psi(e^t)}\int_0^{e^t} f^*(s)\, ds \right)\right|\\
&=\left|-\frac{e^t\psi'(e^t)}{\psi(e^t)}\frac1{\psi(e^t)}\int_0^{e^t} f^*(s)\, ds+ \frac{e^tf^*(e^t)}{\psi(e^t)}\right|\\
&\le 2 \|f\|_{M_\psi}.
\end{align*} 
\end{proof}

The following Theorem strengthens the result from~\cite[Corollary 3.9]{LSS} in the case when $\psi \in \Omega$ satisfies~\eqref{cond_lim}.

\begin{theorem}\label{lemma1}
 Let $\psi \in \Omega$  satisfy~\eqref{cond_lim}. 
A positive operator $T\in \mathcal M_\psi$ is Dixmier measurable if and only if 
there exists a limit in~\eqref{l0}.
\end{theorem}
\begin{proof}
Suppose that $T\in \mathcal M_\psi$ is Dixmier measurable positive operator, 
that is ${\rm Tr}_\omega(T)=A$ for every Dixmier trace ${\rm Tr}_\omega$ on $\mathcal M_\psi$. 
According to Remark~\ref{duality}, we have $\tau_\omega(\mu(T))=A$ for every Dixmier functional $\tau_\omega$ on $M_\psi$. 
Denote, for brevity, $f:=\mu(T)$. 
By Corollary~\ref{add0} we have
\begin{equation}\label{meas}
\lim_{t\to\infty} \frac1{\log t} \int_1^t a(\alpha s,f) \frac{ds}s =A 
\end{equation}
uniformly in $\alpha\ge1$.

Using the pinching theorem one can show that the assumption~\eqref{cond_lim} implies
$$\lim_{t\to\infty} \frac{\psi(Nt)}{\psi(t)}=1 \ \text{for every} \ N>0.$$
So, for any $N>0$ one can find such $t_0=t_0(N)$ that for every $t>t_0$ we have 
\begin{equation}\label{add1}
\frac{\psi(t)}{\psi(Nt)}\ge 1-\frac1N. 
\end{equation}
By the definition of a limit superior, there exists $\alpha>t_0$ such that
\begin{equation}\label{add2}
a(\alpha,f)\ge \left( 1-\frac1N \right) \limsup_{t\to \infty} a(t,f). 
\end{equation}

Using~\eqref{add1} and~\eqref{add2}, we have  
\begin{align*}
 a(s,f)&\ge \frac1{\psi(\alpha N)} \int_0^\alpha f^*(u) du = \frac{\psi(\alpha)}{\psi(\alpha N)} a(\alpha,f)\\
&\ge \left( 1-\frac1N \right)^2 \limsup_{t\to \infty} a(t,f).
\end{align*}
for every $s\in[\alpha, \alpha N]$.

Hence,
$$\frac1{\log N} \int_1^N a(\alpha s,f) \frac{ds}s =
\frac1{\log N} \int_\alpha^{\alpha N} a(s,f) \frac{ds}s \ge \left( 1-\frac1N \right)^2 \limsup_{t\to \infty} a(t,f).$$
Letting $N\to\infty$ and applying~\eqref{meas}, we obtain
$$A \ge \limsup_{t\to \infty} a(t,f).$$

Similarly one can prove that 
$$A \le \liminf_{t\to \infty} a(t,f)$$ 
and, therefore, 
$$\lim_{t\to \infty} a(t,f) =A.$$

The converse implication is trivial.
\end{proof}

Let us consider the Marcinkiewicz space $\mathcal M_\psi$ with $\psi(t)=2^{\sqrt{\log_2 (1+t)}}-1$. 
It is easy to see that $\psi\in \Omega$ satisfies~\eqref{cond_lim}.
Hence, $\mathcal M_\psi$ admits non-trivial Dixmier traces.
A direct computation shows that $\psi(t)=2^{\sqrt{\log_2 (1+t)}}-1$ does not satisfy~\eqref{taub}.

The following Theorem provides an example of a positive operator $T_0\in \mathcal M_\psi$
which is Connes-Dixmier measurable, however it is not Dixmier measurable.

\begin{theorem}\label{4}
Let $\psi(t)=2^{\sqrt{\log_2 (1+t)}}-1$. There exists a positive Connes-Dixmier measurable operator $T_0\in \mathcal M_\psi$ such that
the limit in~\eqref{l0} does not exist.
\end{theorem}
\begin{proof} Let $T_0$ be such that
$$\mu(T_0)=\sup_{k\geq0}2^{k-k^2}\chi_{[0,2^{k^2})}.$$ 
We set $f_0:= \mu(T_0).$ We obtain for every $2^{n^2}\le t<2^{(n+1)^2}$

\begin{equation}\label{a2}
\begin{aligned}
a(t,f_0)&=\frac1{\psi(t)} \left( \int_0^{2^{n^2}} f(s) ds+ (t-2^{n^2}) f_0(t)\right)\\
&=\frac1{\psi(t)} \left( \sum_{k=1}^n  \int_{2^{(k-1)^2}}^{2^{k^2}} 2^{k-k^2}ds+ tf_0^*(t) +O(1)\right) \\
&=\frac{2^{n+1}}{2^{\sqrt{\log_2 t}}} +\frac{tf_0^*(t)}{\psi(t)} +o(1).
\end{aligned} 
\end{equation}

It is easy to see that $a(\cdot,f_0)$ is uniformly bounded and, so, $f_0\in M_\psi$. Hence, $T_0\in \mathcal M_\psi$.

By Lemma~\ref{lem1}, for every dilation invariant generalized limit $\omega$ on $L_\infty$ which is satisfied~\eqref{comp2} we have
\begin{equation}\label{c2}
\omega\left(\frac{tf^*(t)}{\psi(t)}\right)=0 \quad \text{for every} \quad f\in M_\psi. 
\end{equation}

Denote by 
$$x(t):= \sum_{n=0}^\infty 2^{n-\sqrt{\log_2 t}}\chi_{[2^{n^2},2^{(n+1)^2})}(t).$$

We conclude from~\eqref{a2} and~\eqref{c2} that $\tau_\omega(f)= 2\cdot\omega(x)$ 
for every dilation invariant generalized limit $\omega$ on $L_\infty$ satisfying~\eqref{comp2}.

For every $2^{n^2}\le t<2^{(n+1)^2}$, we have
$$(Mx)(t)=\frac1{\log t}\left( \sum_{k=0}^{n-1} \int_{2^{k^2}}^{2^{(k+1)^2}} 2^{k-\sqrt{\log_2 s}}\frac{ds}s +\int_{2^{n^2}}^t 2^{n-\sqrt{\log_2 s}}\frac{ds}s\right).$$
Since 
$$\int 2^{-\sqrt{\log_2 s}}\frac{ds}s=-2\cdot 2^{-\sqrt{\log_2 s}}\left(\sqrt{\log_2 s}+\frac1{\log2}\right)+C$$
and
$$\int_{2^{n^2}}^t 2^{n-\sqrt{\log_2 s}}\frac{ds}s \le \int_{2^{n^2}}^{2^{(n+1)^2}} 2^{n-\sqrt{\log_2 s}}\frac{ds}s = n+ \frac1{\log2} = o(\log t),$$
we have
\begin{align*}
(Mx)(t)&=\frac{-2}{\log t}\left( \sum_{k=0}^{n-1} 2^k \left(2^{-k-1}(k+1+\frac1{\log2})  -2^{-k}(k+\frac1{\log2})\right) \right)+o(1)\\
&=\frac2{\log t} \sum_{k=0}^{n-1} \frac{k}2 +o(1)\\
&=\frac1{2\log2}+o(1).
\end{align*}

Hence, $\lim_{t\to\infty} (Mx)(t)=\frac1{2\log2}$ and, therefore,
$\tau_\omega(f)=\frac1{\log2}$ for every Connes-Dixmier functional $\tau_\omega$.
Consequently, $T_0$ is Connes-Dixmier measurable operator.

However, direct computation shows that
$$\limsup_{t\to\infty} a(t,f_0) \ge \lim_{n\to\infty} a(2^{n^2},f_0)=2$$
 and 
$$\liminf_{t\to\infty} a(t,f_0) \le \lim_{n\to\infty} a(2^{(n+1/2)^2},f_0)=\sqrt{2}.$$
We conclude that $a(\cdot,f_0)$ has no limit at infinity and, so, a limit in~\eqref{l0} does not exist.

\end{proof}

\bibliographystyle{amsplain}

\end{document}